\documentstyle[12pt,leqno]{article}
\begin{document}
\renewcommand{\theequation}{\arabic{section}.\arabic{equation}}
\newcommand{\be}{\begin{eqnarray}}
\newcommand{\en}{\end{eqnarray}}
\newcommand{\no}{\nonumber}
\newcommand{\ti}{\tilde}
\newcommand{\laa}{\lambda}
\newcommand{\la}{\langle}
\newcommand{\ra}{\rangle}
\newcommand{\ga}{\gamma}
\newcommand{\ep}{\epsilon}
\newcommand{\de}{\delta}
\renewcommand{\thefootnote}{}
\newcommand{\pl}{\parallel}
\newcommand{\ov}{\overline}
\newcommand{\bet}{\beta}
\newcommand{\al}{\alpha}
\newcommand{\fr}{\frac}
\newcommand{\pa}{\partial}
\newcommand{\we}{\wedge}
\newcommand{\om}{\Omega}
\newcommand{\na}{\nabla}
\newcommand{\ri}{\rightarrow}
\newcommand{\D}{\Delta}
\newcommand{\pp}{\phi_{\alpha i}}
\newcommand{\ii}{\int_{\Omega}}
\newcommand{\vs}{\vskip0.3cm}
\newcommand{\R}{I\!\!R^{n+1}}
\newcommand{\qed}{\hfill\rule[-1.2mm]{1.2ex}{1.2ex}}

\title {Inequalities for the Steklov Eigenvalues }
\footnotetext{
\hspace*{2ex}Key words and phrases: Steklov eigenvalues, Hadamard manifolds, isoperimetric  bounds, non-negative Ricci curvature, compact manifolds with boundary,  Euclidean ball.}
\author{Changyu Xia\thanks{Partially supported by CNPq and CAPES/PROCAD.} and Qiaoling Wang \thanks{Partially supported by CNPq and CAPES/PROCAD.}  } \date{}
\maketitle
\begin{abstract}
This paper studies eigenvalues of some Steklov problems. Among other things, we show the following sharp estimtes.
Let $\om$ be a  bounded smooth domain in an $n(\geq 2)$-dimensional Hadamard  manifold an let $0=\laa_0 < \laa_1\leq \laa_2\leq ... $ denote the eigenvalues of the Steklov problem: $\D u=0$ in $\om$ and $(\pa u)/(\pa \nu)=\laa u$ on $\pa \om$. Then $\sum_{i=1}^{n} \laa^{-1}_i
\geq (n^2|\om|)/(|\pa\om|) $ with equality holding if and only if $\om$ is isometric to an $n$-dimensional  Euclidean ball.
Let $M$ be an $n(\geq 2)$-dimensional compact connected Riemannian manifold with boundary  and non-negative Ricci curvature.
Assume that the mean curvature of $\pa M$ is bounded below by a positive constant $c$ and let $q_1$ be the first eigenvalue of the Steklov  problem: $ \Delta^2 u= 0$ in $ M$ and $u= (\pa^2 u)/(\pa \nu^2)  -q(\pa u)/(\pa \nu) =0$  on $ \pa M$.
Then  $q_1\geq c$ with equality holding if and only if $M $ is isometric to a ball of radius $1/c$ in ${\bf R}^n$.
\end{abstract}



\section{Introduction}

Let $M^n$ be a compact Riemannian manifold with boundary and let $\rho$ be a positive continuous function on $\pa M$. The second order
Steklov problem is given by
\be\left\{\begin{array}{l}
 \Delta u= 0  \ \ {\rm in \ \ } M, \\
  \fr{\pa }{\pa\nu} u=p\rho u \ \  {\rm on} \ \ \pa M;
  \end{array}\right.
\en
where $\fr{\pa }{\pa\nu}$ is the outward normal derivative and $p$ is a real number.
This problem was first introduced by Steklov in [21] for bounded
domains in the plane where there are several physical interpretations
of the  problem.  In particular, it describes the
vibration of a free membrane with its whole mass $W(M)$ distributed on the
boundary with density $\rho$:
\be\no
W(M)=\int_{\pa M} \rho.
\en
If $\rho\equiv 1$. The mass of $M$ is equal to the area of $\pa M$. The eigenvalues
\be\no
0=\laa_0<\laa_1\leq\laa_2\leq\laa_3\cdots \nearrow +\infty
\en
of the problem (1.1) satisfy the following variational characterization [3]:
\be\no
 \laa_k&=&\min\left\{\fr{\int_{M}|\na v|^2 }{\int_{\pa M} \rho v^2 }: v\in H^{1}(M), \int_{\pa M} \rho v u_j =0, j=1,\cdots, k-1\right\},\\  & & \ \ \ \ \ \ k=1, 2, \cdots
 \en
where $u_0, u_1,\cdots,$ are the eigenfunctions corresponding to the eigenvalues $\laa_0, \laa_1,\cdots,$ respectively. Note that as in the case of Neumann boundary conditions, the Steklov spectrum always starts with the eigenvalue
$\laa_0=0$, and the corresponding eigenfunctions are constant. If $M$ is an $n$-dimensional Euclidean ball of radius $R$ and $\rho \equiv 1$, then $\laa_1$ has multiplicity $n$, the corresponding  eigenfunctions are the coordinate functions $x_i,  i=1,\cdots, n$  and $\laa_1=1/R$.

Many interesting results for the eigenvalues of (1.1) have been obtained during the past years, especially when $M$ is a compact domain in a Euclidean space (cf. [2, 3, 5, 8, 11-20, 24, 25]).
Using conformal mapping, Weinstock [24] has proved the following
isoperimetric inequality: For all two-dimensional simply connected domains with analytic boundary in a Euclidean plane of assigned total
mass $W$, the circle yields the maximum value of $\laa_1$, that is,
\be
\laa_1(M)\leq \fr{2\pi}{\int_{\pa M} \rho },   \ \ {\rm for} \ M\subset {\bf R}^2.
\en
Later on, Hersch and Payne [11] observed that Weinstock's proof actually gives the sharper isoperimetric result
\be
\fr 1{\laa_1(M)}+\fr 1{\laa_2(M)}\geq \fr{\int_{\pa M} \rho }{\pi}, \ \ {\rm for} \ M\subset {\bf R}^2.
\en
In 2000, Brock [5] proved a similar estimate for  compact domains in arbitrary dimensional Euclidean space.

In this paper, we will prove a sharp estimates for the first sum of the reciprocal of the first $n$ non-zero Steklov eigenvalues of a compact domain in a Hadamard manifold. Namely, we have
\vs
{\bf Theorem 1.1}.   {\it Let $Q^n$ be a  Hadamard manifold and let $\om$ be a bounded domain in $Q^n$.
Let $\rho$ be a continuous positive function on $\pa\om$. Then the first $n$ non-zero Steklov eigenvalues of the problem
\be\left\{\begin{array}{l}
\Delta f=0 \  {\rm in \ \ } \om,\\ \left.\fr{\pa f}{\pa \nu}\right|_{\pa \om}=\laa\rho f
\end{array}\right.
\en
satisfy
\be
\sum_{i=1}^n \fr 1{\laa_i}\geq \fr{n^2|\om|}{\int_{\pa\om}  \rho^{-1}},
\en
where $|\om|$ denotes the volume of $\om$. Moreover, when $\rho=\rho_0$ is constant, the equality holding in (1.6) if and only if $\om$ is isometric to an $n$-dimensional Euclidean ball.}
\vs
Consider now the following three fourth order Steklov eigenvalue problems:
\be\left\{\begin{array}{l}
 \Delta^2 u= 0  \ \ {\rm in \ \ } M, \\
u= \D u  -p\fr{\pa }{\pa \nu}u =0 \ \ {\rm on \ \ } \pa M,
\end{array}\right.
\en

\be
\left\{\begin{array}{l}
 \Delta^2 u= 0  \ \ {\rm in \ \ } M, \\
\fr{\pa }{\pa \nu}u = \fr{\pa }{\pa \nu} \D u+\xi u =0 \ \ {\rm on \ \ } \pa M
\end{array}\right.
\en
and
\be\left\{\begin{array}{l}
 \Delta^2 u= 0  \ \ {\rm in \ \ } M, \\
u= \fr{\pa^2 }{\pa \nu^2} u  -q\fr{\pa }{\pa \nu}u =0 \ \ {\rm on \ \ } \pa M.
\end{array}\right.
\en
Elliptic problems with parameters in the boundary conditions are called Steklov problems
from their first appearance in [23]. The problem (1.7) was studied by Kuttler [12] and Payne [18]   who studied the isoperimetric
properties of the first eigenvalue $p_1$.   As pointed out in [12, 13, 15], $p_1$ is the
sharp constant for $L^2$ a priori estimates for solutions of the (second order) Laplace equation
under nonhomogeneous Dirichlet boundary conditions. The whole spectrum
of the biharmonic Steklov problem (1.7) was studied in [9] where one can also find a physical
interpretation of $p_1$. We refer to [4], [6], [25] for some further developments about $p_1$. The problem (1.8) was first studied in [14] where some estimates for the first non-zero eigenvalue $\xi_1$  were proved.  It should be mentioned that the first eigenvalue of the problem  is zero and the corresponding eigenfunctions are constant. The problem (1.9) is a natural Steklov problem and one can check that  when the mean curvature of $\pa M$ is constant, it is  equivalent to (1.7).
 Our second result in this paper is  the following sharp lower bound for the first non-zero eigenvalue of the problem (1.9).
\vs
{\bf Theorem 1.2.} {\it  Let $M^n$ be a compact Riemannian manifold with boundary and non-negative Ricci curvature. Assume that the mean curvature of $\pa M$ is bounded below by a positive constant $c$. Let $q_1$ be the first eigenvalue of the
Steklov  problem (1.9).
Then $q_1\geq c$ with equality holding if and only if $M $ is isometric to a ball of radius $\fr 1c$ in $ {\bf R}^n$.}
\vskip0.3cm
We then prove a lower bound for the first  non-zero eigenvalue of the Steklov problem (1.8).
\vs
{\bf Theorem 1.3.} {\it  Let $M^n$ be a compact connected Riemannian manifold with non-negative Ricci curvature and  boundary. Assume that the principal curvatures of $\pa M$ are bounded from below by a positive constant $c$.  Denote by $\laa_1$ the first non-zero eigenvalue of the Laplacian  acting on functions on $\pa M$.  Then the  first non-zero eigenvalue $\xi_1$ of the
problem (1.8) satisfies $\xi_1 > \fr{nc\laa_1}{n-1} $.}
\vskip0.3cm
Our final result in this paper is concerned with a Steklov eigenvalue problem for elliptic equation in divergence form which generalizes a  result in [23]. Steklov eigenvalue problem for elliptic equations in divergence form on bounded domains in a Euclidean plane has been studied in [1].
\vskip0.3cm
{\bf Theorem 1.4.} {\it  Let $M^n$ be a compact connected Riemannian manifold with boundary $\pa M$. Assume that the Ricci curvature of $M$ is bounded below by $-\kappa_0$ for some non-negative constant $\kappa_0$ and that the principal curvatures of $\pa M$ are bounded  below by a positive constant $c$. Let $A: M\ri {\rm End}(TM)$ be a smooth
symmetric and positive definite section of the bundle of all endomorphisms of $TM$ and assume that there exists a positive number $\delta$ such that  $A\leq \delta I$ in the sense of bilinear form throughout  $M$. Consider the operator $L: C^{\infty}(M)\ri C^{\infty}(M)$
given by $Lu= -{\rm div } (A\na u), \ u\in C^{\infty}(M)$. Denote by $\laa_1$ the first non-zero eigenvalue of the Laplacian  acting on functions on $\pa M$ and let $\eta_1$ be  the first non-zero eigenvalue  of the following
Steklov eigenvalue problem :
\be
\left\{\begin{array}{l}
  Lf=0 \ \ \ \ \ {\rm in}\ \ \  M, \\
\fr{\pa }{\pa \nu}u -\eta u =0 \ \ {\rm on \ \ } \pa M
\end{array}\right.
\en
Then we have
 \be
(2\laa_1+\kappa_0)^2\geq 4(n-1)\laa_1c^2
\en
and
\be
\eta_1\leq \fr{2\laa_1+\kappa_0+\sqrt{(2\laa_1+\kappa_0)^2-4(n-1)\laa_1c^2}}{2(n-1)c}\delta.
\en
Furthermore, if (1.11) or (1.12) take equality sign, then $M$ is isometric  to an $n$-dimensional Euclidean ball of radius $\fr 1c$.}
\vskip0.3cm

\section{Proof of the Results}
\setcounter{equation}{0}
In this section, we will prove the main results as mentioned in the last section.
\vskip0.3cm

{\it Proof of Theorem 1.1.}
We denote by $\la, \ra$  the Riemannian metric on $M$. For any $p\in M$, let ${\rm exp}_p$ and
$UM_p$ be the exponential map and unit tangent space of $M$ at $p$, respectively. Let
$\{u_i\}_{i=0}^{+\infty}$  be an orthonormal set of eigenfunctions corresponding to the eigenvalues $\{\laa_i\}_{i=0}^{+\infty}$
 of the problem (1.1). For $  i=1, 2,\cdots,n,$
 we need to choose trial functions $\phi_i$ for each of the eigenvalues $\laa_i$ and
insure that these are orthogonal to the  eigenfunctions $u_0, u_1, \cdots, u_{i-1}$.
To do this, we let $p\in M, \{e_1,\cdots, e_n\}$  be an orthonormal basis of $T_pM$, the tangent
space of $M$ at $p$, and $y : M\ri {\bf R}^n$ be the Riemannian normal coordinates on $M$
determined by $(p; e_1,\cdots,e_n)$. It follows from  the Cartan-
Hadamard theorem [7] that $y$ is defined on all of $M$ and is a
diffeomorphism of $M$ onto ${\bf R}^n$. We can choose $p$ and $\{e_1,\cdots, e_n\}$ so that the respective coordinate
functions $ y^i : M\ri {\bf R}, i=1,\cdots,n,$ of $y : M \ri {\bf R}^n$ satisfy
\be
\int_{\pa \om} \rho y^i =0.
\en
In fact, parallel translate the frame $\{e_1,\cdots, e_n\}$ along every geodesic emanating
from $p$ and thereby obtain a differentiable orthonormal frame field $\{E_1, \cdots, E_n\}$ on $M$. Let $y_q : M\ri {\bf R}^n$
denote the Riemann normal coordinates of $M$ determined by $\{E_1, \cdots, E_n\}$ at $q$, and let $y_q^i, i=1,\cdots, n,$ be the coordinate functions of $y_q$. By
definition, for $i=1,\cdots,n,  y_q^i : M\ri {\bf R}$ is given by $y^i_q(z)=\la {\rm exp}_q^{-1}(z), E_i(q)\ra$, then
\be
Y(q)=\sum_{i=1}^n \left\{\int_{\pa\om}\rho y_q^i \right\} E_i(q)
\en
is a continuous vector field on $M$. If we restrict $Y$ to a geodesic ball $B$ containing $\overline{\om}$
then the convexity of $B$ implies that on the boundary of $B$; $Y$ points into $B$. The
Brouwer fixed point theorem  then implies that $X$ has a zero. So we may assume
that $p$ and $\{e_1,\cdots, e_n\}$,  actually satisfy (2.1).

For any $w\in UM_p$, let $\theta_w$ be the function on $M$ defined by $\theta_w(z)=\la {\rm exp}_p^{-1}(z), w\ra $. Then (2.1) is equivalent to say that for $i=1, \cdots, n,$
\be
\int_{\pa\om} \rho\theta_{e_i} =0.
\en
Thus for any $\sigma\in UM_p$,
\be
\int_{\pa\om} \rho\theta_{\sigma}=0.
\en
Next we argue, via the Borsuk-Ulam theorem [22, p. 266] that there
exist $(n-1)$ unit orthogonal vectors $\sigma_2,\cdots,\sigma_n$ in $T_pM$ such that
\be
\int_{\pa\om} \rho\theta_{\sigma} u_i=0,
\en
for $j=2, 3,\cdots,n$ and $i=1,\cdots, j-1.$ To see this, we define a mapping
$g_n: UM_p\ri {\bf R}^{n-1}$ componentwise by
\be
g_{n,k}: \sigma\ri \int_{\pa\om}\rho\theta_{\sigma_j} u_k, \ \ {\rm for}\ k=1,\cdots, n-1.
\en
Since  $UM_p$ is isometric to an $(n-1)$-dimensional Euclidean sphere and $g_n(-\sigma)=-g_n(\sigma), \ \forall \sigma\in UM_p$, the
Borsuk-Ulam theorem tells us that there is a direction
$\overline{\sigma}\in UM_p$ such that $g_n(\overline{\sigma})=0.$ We take this $\overline{\sigma}$ as our $\sigma_n$. So we have found a
$\sigma_n\in UM_p$  such that
\be
\int_{\pa\om}\rho\theta_{\sigma_n} u_i =0, \  i=1,\cdots, n-1.
\en
Let $S_{\sigma_n}^{n-2}$ be the equator determined by $\sigma_n$ in $UM_p$, that is, $S_{\sigma_n}^{n-2}=
\{\sigma\in UM_p|\la \sigma, \sigma_n\ra =0\}$; then $S_{\sigma_n}^{n-2}$ is isometric to an $(n-2)$-dimensional unit
Euclidean sphere. By the same arguments as above, we can find  a $\sigma_{n-1}$ in $S_{\sigma_n}^{n-2}$ such that
\be
\int_{\pa\om}\rho\theta_{\sigma_{n-1}} u_i =0, \  i=1,\cdots, n-2.
\en
Continuing in this way, we get finally $(n-1)$ mutually orthogonal unit tangent vectors  $\sigma_2,\cdots,\sigma_n\in UM_p$ satisfying (2.5).
We extend $\{\sigma_2,\cdots,\sigma_n\}$ to be an orthonormal basis of $T_pM : \{\sigma_1, \sigma_2,\cdots,\sigma_n\}$.
Let $x : M\ri {\bf R}^n$ be the Riemannian normal coordinates on $M$ determined by $(p; \sigma_1,\cdots,\sigma_n);$ then the respective coordinate functions $x_i: M\ri {\bf R}, i=1,\cdots, n$ of $x : M\ri {\bf R}^n$ are given by $x_i(q)=\la{\rm exp}_p^{-1}(q), \sigma_i\ra$. Since $u_0$ is just
the constant function $(\int_{\pa\om}\rho)^{-1/2}$, we  have from (2.4) and (2.5) that
\be
\int_{\pa \om} \rho x_j u_i =0,
\en
for $j=1,\cdots, n$ and $i=0, 1,\cdots, j-1.$

For the $n$ trial functions $\phi_1,\cdots, \phi_n$ we simply   choose
\be
\phi_j= x_j \ \ {\rm for} \ j=1,\cdots,n.
\en
It then follows from (1.2) and (2.9) that
\be
\laa_i\leq \fr{\int_{ \om} |\na x_i|^2 }{\int_{\pa \om} \rho x_i^2 }, \ \ i=1,\cdots, n.
\en
Let $\left\{\fr{\pa}{\pa x_k}, k=1,\cdots,  n\right\}$ be the natural basis of the tangent spaces associated with
the coordinate chart $x$ and let $g_{kl}=\left\la \fr{\pa}{\pa x_k}, \fr{\pa}{\pa x_l}\right\ra $, $k, l=1,\cdots, n$.
Since $M$ has
non-positive sectional curvature, the Rauch comparison theorem [7] implies that
the eigenvalues of the matrix $(g_{kl})$ are all $\geq 1$. Thus the eigenvalues of
$(g^{kl})=: (g_{kl})^{-1}$ are $\leq 1$ and so we have $g^{kk}\leq 1, \ k=1,\cdots, n$.
Thus we have for $i=1,\cdots, n,$
\be
|\na x_i|^2=\left\la \sum_{k=1}^n g^{ik} \fr{\pa}{\pa x_k}, \sum_{l=1}^n g^{il} \fr{\pa}{\pa x_l}\right\ra= g^{ii} \leq 1.
\en
Substituting (2.12) into (2.11), we get
\be
\laa_i\leq \fr{|\om|}{\int_{\pa\om}\rho x_i^2 }, \ i=1,\cdots, n.
\en
By inverting and summing on $i$,
\be
\sum_{i=1}^n \fr 1{\laa_i}&\geq& \fr{\int_{\pa\om}\rho\sum_{i=1}^n x_i^2 }{|\om|}=
\fr{\int_{\pa \om} \rho r^2 }{|\om|}
\en
where $r = d(p, \cdot): M\ri {\bf R}$ denotes the distance function from $p$. The Cauchy-Schwart inequality implies that
\be
\int_{\pa\om} \rho r^2 \geq\fr{(\int_{\pa \om} r )^2}{\int_{\pa\om} \rho^{-1}}
\en
with equality holding if and only if $r\rho=const.$ on $\pa\om$. On the other hand, since $M$ is a Hadamard manifold, the Laplace comparison theorem
tells us that (cf. [21])
\be
\D r^2 \geq 2n.
\en
Integrating (2.16) on $\om$ and using the divergence theorem, we get
\be n|\om|\leq \int_{\pa \om} r\la \na r, \nu\ra\leq \int_{\pa \om} r |\na r|=\int_{\pa \om} r.
\en
Combining (2.14), (2.15) and (2.17), we get (1.6).

Assume now that $\rho=\rho_0$ is constant and that  the equality holds in (1.6). In this case, we must have $r|_{\pa\om}=const.$
and so $\om$ is a geodesic ball with center $p$. Also, we have
\be
\D r^2|_{\om}= 2n.
\en
It then follows from the equality case in the Laplace comparison theorem [21] and the Cartan's theorem [7] that $\om$ is isometric to an $n$-dimensional
Euclidean ball. This completes the proof of Theorem 1.1. \qed\\
\vskip0.2cm
Before proving theorems 1.2-1.4, we first recall Reilly's formula.
Let $M$ be an $n$-dimensional compact manifold  with boundary. We will often write $\la, \ra$ the Riemannian metric on $M$ as well as that induced on $\pa M$. Let $\na$ and $\D $ be the connection  and the Laplacian on $M$,
respectively. Let $\nu$ be the unit outward normal vector of $\pa M$. The shape operator of $\pa M$ is given by $S(X)=\na_X \nu$ and the second fundamental form of $\pa M$ is defined as $II(X, Y)=\la S(X), Y\ra$, here $X, Y\in T \pa M$. The eigenvalues of $S$ are called the principal curvatures of $\pa M$ and the mean curvature $H$ of $\pa M$ is given by $H=\fr 1{n-1} {\rm tr\ } S$, here ${\rm tr\ } S$ denotes the trace of $S$.  We can now state the Reilly's formula  (see [20, p. 46]).
For a smooth function $f$ defined on  $M$, the following identity holds if
$h=\left.\fr{\pa }{\pa \nu}f\right|_{\pa M}$, $z=f|_{\pa M}$ and ${\rm Ric}$ denotes the Ricci tencor of $M$:
\be
& & \int_M \left((\D f)^2-|\na^2 f|^2-{\rm Ric}(\na f, \na f)\right)\\ \no
&=& \int_{\pa M}\left( ((n-1)Hh+2\overline{\D}z)h + II(\overline{\na}z, \overline{\na}z)\right).
\en
Here $\na^2 f$ is the Hessian of $f$; $\ov{\D}$ and $\ov{\na}$  represent the Laplacian and the gradient on $\pa M$ with respect
to the induced metric on $\pa M$,  respectively.
\vskip0.3cm
{\it Proof of Theorem 1.2.}
Let $w$ be an eigenfunction corresponding to the first eigenvalue $q_1$  of the problem (1.9):
\be\left\{\begin{array}{l}
 \Delta^2 w= 0  \ \ {\rm in \ \ } M, \\
w= \fr{\pa^2 }{\pa \nu^2}w  -q_1\fr{\pa }{\pa \nu}w =0 \ \ {\rm on \ \ } \pa M.
\end{array}\right.
\en
Observe that $w$ is not a constant. Otherwise, we would conclude from $w|_{\pa M}=0$ that $w\equiv 0$. Set $\eta=\fr{\pa }{\pa \nu}w|_{\pa M}$;
then $\eta\neq 0$. In fact, if $\eta=0$, then
$$w|_{\pa M}=\na w|_{\pa M}=\left.\fr{\pa^2 }{\pa \nu^2}w\right|_{\pa M}=0$$
which implies that $\Delta w|_{\pa M}=0$ and so $\Delta w=0$ on $M$ by the maximum principal [10], which in turn implies that $w=0$. This is a contradiction.

Since $w|_{\pa M}=0$, we obtain from the divergence theorem that
\be
0=\int_M \langle \na w, \na(\Delta w)\rangle=-\int_M w\Delta^2 w=0.
\en
It follows that
\be
\int_{\pa M} \D w\fr{\pa }{\pa \nu}w=\int_M\langle\na (\D w), \na w\rangle +\int_M (\D w)^2=\int_M (\D w)^2.
\en
Since
\be
\D w|_{\pa M}&=&\fr{\pa^2 }{\pa \nu^2}w+(n-1)H\fr{\pa }{\pa\nu}w\\ \no
&=& q_1\fr{\pa }{\pa\nu}w+(n-1)H\fr{\pa }{\pa\nu}w,
\en
we conclude that
\be
q_1=\fr{\int_M(\D w)^2-(n-1)\int_{\pa M} H\eta^2}{\int_{\pa M} \eta^2}
\en
Substituting $w$ into  Reilly's formula, we have
\be
\int_{M} \left\{ (\D w)^2- |\na^2 w|^2\right\}&= &\int_{M}
{\rm Ric}(\na w, \na w)+\int_{\pa M} (n-1)H \eta^2\\ \no
&\geq & (n-1)\int_{\pa M}H \eta^2.
\en
The Schwarz inequality implies that
\be
|\na^2 w|^2\geq \fr 1n(\D w)^2
\en
with equality holding if and only if $\na^2 w=\fr{\D w}n \la, \ra$.
Substituting (2.26) into (2.25), we have
\be
\int_M (\D w)^2 \geq n\int_{\pa M} H \eta^2
\en
which, combining with (2.24), gives
\be
q_1\geq \fr{\int_{\pa M} H\eta^2}{\int_{\pa M} \eta^2}\geq c.
\en
Assume now that   $q_1 = c$. In this case, the inequalities (2.25) and (2.28)  must take equality sign. In particular, we have
\be
 \ \ \na^2 w=\fr {\D w}{n}\langle, \rangle.
\en
Take an orthornormal frame $\{e_1,\cdots, e_{n-1}, e_{n}\}$ on $M$ such that when restricted to $\pa M$, $e_{n}=\nu$.
 From
$0=\na^2 w(e_i, e_{n})$, $i=1, \cdots, n-1$, and $w|_{\pa M}=0$, we conclude that
$\eta=b_0=const.$  Since (2.28) takes equality sign and $\eta $ is constant, we infer that $H\equiv c$
which implies from (2.23) that $\D w|_{\pa M}=ncb_0$. Consequently, we conclude from the maximum principle that
$\D w$ is constant on $M$.   Without loss of generality, let us assume that $\D w=1$ and so
\be
 \ \ \na^2 w=\fr 1n\langle, \rangle.
\en
Deriving (2.30) covariantly, we get
 $\na^3 w=0$ and from the Ricci identity,
\be
R(X, Y)\na w=0,
\en
for any tangent vectors  $X, Y$ on $M$, where $R$ is the curvature tensor of $M$. By the the maximum principle $w$ attains its minimum at some point $x_0$ in the interior of $M$. Let $r$ be the distance function to $x_0$; then by  (2.30) it follows that
\be
\na w=\fr 1n r \fr{\pa}{\pa r}.
\en
Using (2.31), (2.32), Cartan's theorem (cf. [7]) and $w|_{\pa M}=0$, we conclude that $M$ is an Euclidean ball whose center is $x_0$, and
\be\no
w(x)=\fr 1{2n}(|x-x_0|^2-b^2)
\en
in $M$, here $b$ is the radius of the ball. Since the mean curvature of $\pa M$ is $c$,  the radius of the ball is $\fr 1c$. This completes the proof of Theorem 1.2.\qed\\
\vs
{\it Proof of Theorem 1.3.}
Let $f$ be an eigenfunction corresponding to the first non-zero eigenvalue $\xi_1$ of the problem (1.1), that is
\be
\left\{\begin{array}{l} \Delta^2 f= 0  \ \ {\rm in \ \ } M, \\
\fr{\pa }{\pa \nu} f = \fr{\pa }{\pa \nu} \D f+\xi_1 f =0 \ \ {\rm on \ \ } \pa M.
\end{array}
\right.
\en
Set $z=f|_{\pa M}$; then $z\neq 0$ and
\be
\xi_1=\fr{\int_{M}(\D f)^2}{\int_{\pa M} z^2}.
\en
Substituting $f$ into  Reilly's formula, we have
\be
\int_{M} \left\{ (\D f)^2- |\na^2 f|^2\right\}&=& \int_M {\rm Ric}(\na f, \na f)+  \int_{\pa M} II(\overline{\na}z, \ov{\na} z)
\\ \no &\geq & c\int_{\pa M} |\ov{\na}z|^2.
\en
The Schwarz inequality implies that
\be
|\na^2 f|^2\geq \fr 1n(\D f)^2
\en
with equality holding if and only if $\na^2 f=\fr{\D f}n \la, \ra$.

From
\be
\D^2 f=0\ \ \ \  {\rm on}\ \ M,
\en
we have
\be
\int_{\pa M}\fr{\pa }{\pa \nu}\D f=0
\en
and so we conclude from (2.33) that $\int_{\pa M } z=0$. It then follows from the Poincar\'e inequality that
\be
\int_{\pa M} |\na z|^2\geq \laa_1 \int_{\pa M}  z^2.
\en
Combining (2.34)-(2.36), (2.39), we get $\xi_1\geq \fr{nc\laa_1}{n-1}. $
Let us show by contradiction that the case  $\xi_1= \fr{nc\laa_1}{n-1} $
can't occur. In fact, if   $\xi_1= \fr{nc\laa_1}{n-1} $, then
we must have
\be
\na^2 f=\fr{\D f}n \la, \ra
\en
If $X$ is a tangent vector field along $\pa M$, then we have from (2.40) and $\left.\fr{\pa }{\pa \nu}f\right|_{\pa M}=0$ that
\be
0=\na^2f(\nu, X)=X\nu f-(\na_X\nu)f=-\la \na_X \nu, \ov{\na} z\ra.
\en
In particular, we have
\be\no
II(\ov{\na }z, \ov{\na }z)=0.
\en
This is impossible since $II\geq c I$ and $ z$ is not  constant. The proof of Theorem 1.3 is completed.\qed \\
\vs

{\it Proof of Theorem 1.4.}
Let $g$ be the solution of the following Laplace equation
\be\no
 \left\{\begin{array}{l}
  \D g=0 \ \ \ \ \ \mbox{in}\ \ \  M, \\
  g|_{\pa M}=v,
\end{array}\right.
\en
where $v$ is a first eigenfunction of $\pa M$ corresponding to $\laa_1$, that is $\ov{\D} v+\laa_1 v=0$. Set $h=\left.\fr{\pa }{\pa \nu}g\right|_{\pa M}$.
Sine the principal curvature of $\pa M$ are bounded below by $c$, we have
\be
II(\ov{\na}v, \ov{\na}v)\geq c|\ov{\na} v|^2, \ \ H\geq c.
\en
It then follows by
substituting $g$ into the Reilly's formula and noticing that the Ricci curvature of $M$ is bounded below by $-\kappa_0$ that
\be& &
\kappa_0\int_M   |\na g|^2\\ \no &\geq& \int_M \left((\D g)^2-|\na^2 g|^2-{\rm Ric}(\na g, \na g)\right)\\ \no
&\geq & (n-1)c\int_{\pa M} h^2 -2\laa_1 \int_{\pa M} hv +c \int_{\pa M}|\ov{\na} v|^2\\ \no
&= & (n-1)c\int_{\pa M} h^2 -2\laa_1 \int_{\pa M} hv +c\laa_1 \int_{\pa M} v^2.
\en
Since
\be
\int_M   |\na g|^2=\int_{\pa M} hv,
\en
we have
\be
0 \geq (n-1)c\int_{\pa M} h^2 -(2\laa_1+\kappa_0) \int_{\pa M} hv +c\laa_1 \int_{\pa M} v^2
\en
which gives
\be
0 &\geq& (n-1)c\int_{\pa M} \left(h -\fr{(2\laa_1+\kappa_0)}{2(n-1)c}v\right)^2\\ \no & &\ \ \  +\left(c\laa_1-\fr{(2\laa_1+\kappa_0)^2}{4(n-1)c}\right) \int_{\pa M} v^2\\ \no &\geq& \left(c\laa_1-\fr{(2\laa_1+\kappa_0)^2}{4(n-1)c}\right) \int_{\pa M} v^2.
\en
Thus we have
$$
c\laa_1-\fr{(2\laa_1+\kappa_0)^2}{4(n-1)c}\leq 0.
$$
This proves (1.11).
\vs
Observe that  $v=g|_{\pa M}$ is an eigenfunction of the Laplacian of $\pa M$, and so we have $\int_{\pa M} v=0$. It then follows from the variational characterization of $\eta_1$ (cf. [1]) that
\be
\eta_1&\leq &\fr{\int_M \langle A\na g, \na g\rangle}{\int_{\pa M} v^2}\\ \no
&\leq &\fr{\delta \int_M |\na g|^2}{\int_{\pa M} v^2}
\\ \no
&=&\fr{\delta \int_{\pa M} vh}{\int_{\pa M} v^2}
\\ \no
&\leq &\delta\left\{\fr{\int_{\pa M} h^2}{\int_{\pa M} v^2}\right\}^{1/2}
\en
It follows from (2.45) that
\be\no
0 \geq (n-1)c\int_{\pa M} h^2 -(2\laa_1+\kappa_0) \left(\int_{\pa M} h^2\right)^{1/2}\left(\int_{\pa M} v^2\right)^{1/2} +c\laa_1 \int_{\pa M} v^2
\en
and so
\be& &
\left(\int_{\pa M} h^2\right)^{\fr 12}\\ \no &\leq& \fr{2\laa_1+\kappa_0+\sqrt{(2\laa_1+\kappa_0)^2-4(n-1)\laa_1c^2}}{2(n-1)c}\left(\int_{\pa M} v^2\right)^{\fr 12}.
\en
Substituting (2.48)   into (2.47), one gets (1.12).

This finishes the proof of the first part of Theorem 1.1. If (1.10) take equality sign, then  the inequalities in (2.43) and (2.46) should also   take equality sign. Thus we have
\be
\na^2 g=0, \ \  \ \ (H-c)h=0
\en
and
\be
h=\fr{2\laa_1+\kappa_0}{2(n-1)c}v.
\en
Take a local orthonormal fields $\{e_i\}_{i=1}^{n-1}$ tangent to $\pa M$. We know from (2.49) and (2.50) that
\be\no
0&=&\sum_{i=1}^{n-1}\na^2 g(e_i, e_i)=\ti{\D} z +(n-1)H h
\\ \no
&=& -\laa_1 v+(n-1)c\cdot\fr{2\laa_1+\kappa_0}{2(n-1)c}v,
\en
which gives $\kappa_0=0$ and $\laa_1=(n-1)c^2$. It then follows from a result in [26]  that $M $ is isometric to an $n$-dimensional Euclidean ball of radius
$\fr 1c$.

Similarly, if (1.12) take equality sign, then (2.49) should hold and we have
\be
h=\fr{2\laa_1+\kappa_0+\sqrt{(2\laa_1+\kappa_0)^2-4(n-1)\laa_1c^2}}{2(n-1)c}v
\en
which gives again that $\kappa_0=0$ and $\laa_1=(n-1)c^2$. Consequently,  $M $ is isometric to an $n$-dimensional Euclidean ball of radius
$\fr 1c$. This completes the proof of Theorem 1.4.
\vs

\vskip0.4cm
\noindent  Changyu Xia ( xia@mat.unb.br ) and Qiaoling Wang ( wang@mat.unb.br )

\noindent Departamento de Matem\'atica, Universidade de Bras\'{\i}lia, 70910-900 Bras\'{\i}lia-DF, Brazil

\end{document}